# Linear Optimization for the Perfect Meal: A Data-Driven Approach to Optimising the Perfect Meal Using Gurobi

Utkarsh Prajapati[#1], Tanushree Jain[#2], Abhishek Machiraju[#3], Divyam Kaushik[#4]

[#]*IEOR Department, Columbia University*
*USA*

[1]up2147@columbia.edu
[2]tj2555@columbia.edu
[3]rcm2194@columbia.edu
[4]dk3316@columbia.edu

*Abstract*— **This study aims to optimize meal planning for nutritional health and cost efficiency using linear programming. Linear optimization provides an effective framework for addressing the problem of an optimal diet, as the composition of food can be naturally modeled as a linearly additive system. Leveraging a comprehensive nutrition dataset, our model minimizes meal costs while meeting specific nutritional requirements. We explore additional complexities, such as fractional weights and nutrient ratio constraints, enhancing the robustness of the solution. Case studies address common nutritional challenges, providing tailored diet plans. The significance lies in aiding individuals to form balanced, cost-effective dietary schedules, considering fitness goals and caloric needs. This research contributes to efficient, sustainable, and time-sensitive meal planning, emphasizing the intersection of nutrition, optimization, and real-world applicability.**

## 1. Introduction

Eating a healthy diet is crucial for many reasons. general health and the avoidance of noncommunicable diseases[3], improved athletic & cognitive performance[2] and more. Moreover, meals that are well-planned have also been found to be extremely important. For example, in the workplace, they have been shown to increase worker morale and productivity, decrease worker accidents, and, if made economically feasible (cost-efficient), even prove to be a successful company strategy.[1].

In our quest to revolutionize meal planning for optimal nutrition and cost efficiency, our study employs a robust optimization approach[5]. We explore dietary optimization by utilizing the capabilities of linear programming, since nutrition very naturally emerges as a problem that can be modeled as nutrients being ingested in a person in a linearly additive system. The basis of the study is an extensive nutrition dataset that has been painstakingly processed, guaranteeing precision in our estimates. Our approach improves upon conventional meal planning by adding limits on nutrient ratios and fractional weights, which increases the accuracy of our answers [6].

Our study goes beyond the technical details to practical applications, using case studies to address common nutritional challenges. These case studies examine common problems such as vitamin D deficiency.[4], iron deficiency[7], calorie deficit diets[8], and more, emphasizing the practical implications of our research.

Our study aims to provide customized solutions that meet individual needs, rather than just diet plans. We want to enable people to create healthy, economical eating plans by taking into account things like their fitness objectives and energy needs. The confluence of nutrition, optimisation, and practical application drives our research in the direction of effective, long-lasting, and time-efficient meal planning.

## 2. Background

### 2.1 Dataset

The dataset we utilized is the Nutrition dataset [9]. Some background on this dataset - it contains a variety of ingredients required for making a diet - The data points are of wide variety with the broadest classifications being:

Raw food - containing the likes of meat, seafood, vegetables, fruit etc.

Cooked food - containing the likes of scrambled egg, cooked meat and seafood(fried fish, with finer distinction of dry heat).

Beverages - containing the likes of soups, coffee, tea, soft beverages, milkshakes, smoothies, alcoholic beverages.

Prepared food items from popular restaurant chains - like Applebee's crunchy onion rings, Burger King french toast sticks etc.

Common branded snacks - Keebler Danish Wedding Cookies, Goya Crackers etc.

Each data point has the following(72) nutritional components attached to it with their corresponding units for a serving size of (originally) 100g:

| | | | |
|---|---|---|---|
| alanine (g)', | fiber (g)', | manganese (mg)', | sugars (g)', |
| alcohol (g)', | folate (mcg)', | methionine (g)', | theobromine (mg)' |
| arginine (g)', | folic_acid (mcg)', | monounsaturated_fatty_acids (g)', | thiamin (mg)', |
| ash (g)', | fructose (g)', | niacin (mg)', | threonine (g)', |
| aspartic_acid (g)', | galactose (g)', | pantothenic_acid (mg)', | tocopherol_alpha (mg)', |
| caffeine (mg)', | glucose (g)', | phenylalanine (g)', | total_fat (g)', |
| calcium (mg)', | glutamic_acid (g)', | phosphorous', | tryptophan (g)', |
| calories (kcal)', | glycine (g)', | polyunsaturated_fatty_acids (g)', | tyrosine (g)', |
| carbohydrate (g)', | histidine (g)', | potassium (mg)', | valine (g)', |

| carotene_alpha (mcg)', | hydroxyproline (g)', | proline (g)', | vitamin_a (IU)', |
| --- | --- | --- | --- |
| carotene_beta (mcg)', | iron (mg)', | protein (g)', | vitamin_a_rae (mcg)', |
| cholesterol (mg)', | isoleucine (g)', | riboflavin (mg)', | vitamin_b12 (mg)', |
| choline (mg)', | lactose (g)', | saturated_fat (g)', | vitamin_b6 (mg)', |
| copper (mg)', | leucine (g)', | saturated_fatty_acids (g)', | vitamin_c (mg)', |
| cryptoxanthin_beta (mcg)', | lutein_zeaxanthin (mcg)', | selenium (mcg)', | vitamin_d (IU)', |
| cystine (g)', | lysine (g)', | serine (g)', | vitamin_e (mg)', |
| fat (g)', | magnesium (mg)', | sodium (mg)', | vitamin_k (mcg)', |
| fatty_acids_total_trans (mg)', | maltose (g)', | sucrose (g)', | zinc (mg)', |

Figure 1: Overview of dataset used

Typically head of the data frame looks as follows:

| | name | serving_size (g) | calories (kcal) | total_fat (g) | saturated_fat (g) | cholesterol (mg) | sodium (mg) | choline (mg) | folate (mcg) | folic_acid (mcg) | ... | fat (g) | saturated_fatty_acids (g) | monounsaturate |
| --- | --- | --- | --- | --- | --- | --- | --- | --- | --- | --- | --- | --- | --- | --- |
| 0 | Cornstarch | 1.0 | 3.81 | 0.001 | 0.000 | 0.00 | 0.09 | 0.004 | 0.00 | 0.0 | ... | 0.0005 | 0.00009 | |
| 1 | Nuts, pecans | 1.0 | 6.91 | 0.720 | 0.062 | 0.00 | 0.00 | 0.405 | 0.22 | 0.0 | ... | 0.7197 | 0.06180 | |
| 2 | Eggplant, raw | 1.0 | 0.25 | 0.002 | 0.000 | 0.00 | 0.02 | 0.069 | 0.22 | 0.0 | ... | 0.0018 | 0.00034 | |
| 3 | Teff, uncooked | 1.0 | 3.67 | 0.024 | 0.004 | 0.00 | 0.12 | 0.131 | 0.00 | 0.0 | ... | 0.0238 | 0.00449 | |
| 4 | Sherbet, orange | 1.0 | 1.44 | 0.020 | 0.012 | 0.01 | 0.46 | 0.077 | 0.04 | 0.0 | ... | 0.0200 | 0.01160 | |

--

| g) | saturated_fatty_acids (g) | monounsaturated_fatty_acids (g) | polyunsaturated_fatty_acids (g) | fatty_acids_total_trans (mg) | alcohol (g) | ash (g) | caffeine (mg) | theobromine (mg) | price_per_unit |
| --- | --- | --- | --- | --- | --- | --- | --- | --- | --- |
| 5 | 0.00009 | 0.00016 | 0.00025 | 0.00 | 0.0 | 0.0009 | 0.0 | 0.0 | 0.002181 |
| 7 | 0.06180 | 0.40801 | 0.21614 | 0.00 | 0.0 | 0.0149 | 0.0 | 0.0 | 0.003998 |
| 8 | 0.00034 | 0.00016 | 0.00076 | 0.00 | 0.0 | 0.0066 | 0.0 | 0.0 | 0.000283 |
| 8 | 0.00449 | 0.00589 | 0.01071 | 0.00 | 0.0 | 0.0237 | 0.0 | 0.0 | 0.003995 |
| 0 | 0.01160 | 0.00530 | 0.00080 | 0.01 | 0.0 | 0.0040 | 0.0 | 0.0 | 0.000894 |

Figure 2: Dataframe head on Jupyter Notebook

Below are the examples of some of the rather piquant data points and the corresponding 72 nutritional features contained in them:

```
name                    McDONALD'S, BIG MAC
serving_size (g)                        1.0
calories (kcal)                        2.57
total_fat (g)                          0.15
saturated_fat (g)                     0.038
                          ...
alcohol (g)                             0.0
ash (g)                              0.0187
caffeine (mg)                           0.0
theobromine (mg)                        0.0
price_per_unit                     0.004992
Name: 240, Length: 75, dtype: object
```

Figure 3: Datapoint - A McDonald's Big Mac

```
name                  KEEBLER, Iced Oatmeal Cookies
serving_size (g)                                 1.0
calories (kcal)                                 4.67
total_fat (g)                                   0.18
saturated_fat (g)                              0.058
                             ...
alcohol (g)                                      0.0
ash (g)                                          0.0
caffeine (mg)                                    0.0
theobromine (mg)                                 0.0
price_per_unit                              0.004513
Name: 736, Length: 75, dtype: object
```

Figure 4 : Datapoint - One serving of Keebler Iced Oatmeal Cookies

**2.2 Nutritional Deficiencies**

Nutritional shortages arise when the body does not have enough vitamins, minerals, or other nutrients to function properly. Iron deficiency causes anemia and fatigue; vitamin D deficiency affects bone health and immunity; vitamin B12 deficiency causes exhaustion and neurological difficulties; and iodine deficiency causes thyroid abnormalities. Inadequate calcium consumption can have an effect on bone strength, whereas insufficient folate can damage cell division. Modern diets are frequently deficient in omega-3 fatty acids, which are essential for heart health and cognitive function. It is critical to address these deficiencies through balanced meals, supplements, or fortified foods in order to avoid health issues and maintain overall well-being. Regular health screenings can assist in identifying and addressing dietary deficiencies [19].In our study we have outlined some cases of deficiencies and how the diet can be catered accordingly.

# 3. Methods

**3.1 Data Scraping:**

**3.1.1 Assumptions for scraping through Instacart:**
In order for us to solve the problem of minimizing the cost of the meal, we required the price of each and every one of the data points that was given in the dataset which we had to scrape from the internet. For scraping the data, first of all, we took the following assumptions:

1. Results have been scraped from data present at the Instacart app 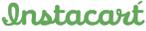 .

2. The prices have been normalized to the most common representation of per unit of the ingredients. For example, prices for apples would be per one apple, whereas prices for sugar would be per gram.

3. The above further assumes that the minimum amount of any ingredient which can be used without wasting food would be a unit.

4. Similar ingredients have similar prices. The data was scraped making use of the 1st two elements of the ingredient name. Since we had a huge number of ingredients, we clustered similar ingredients under a similar first name and scraped the price for a product similar to that first name. For example, corn salted, and corn not-salted would have the same price in our case, also since the difference between both of the ingredients is pretty less - in price as well as in nutritional value.

5. The search results from instacart are assumed to have the best possible search engine optimization. The mapping of the ingredient name to the product whose price has been taken is assumed to be taken care of by Instacart. So we scrape the price of the 1st product which shows up after searching on instacart.

### 3.1.2 Scraping via Selenium:
Utilizing Selenium's Webdriver library[20]. As explained in the assumptions above, the scraping is done from Instacart assuming that it would be the best source for finding all or most of the ingredients in our data. We make use of the search engine optimization algorithms Instacart uses to get the best estimate of the price per unit for each of the ingredients.

But before searching for the ingredients on Instacart, the name of the ingredients had to be cleaned up to reduce the complexity of the model. Similar ingredients were clustered together and were assumed to have similar price as well as nutritional value. Similar ingredients were defined based on the words in the name of the ingredients. Ingredients having "corn" in their name would be clubbed together.

Step 1: Selenium access Instacart through the webdriver.
Step 2: For each ingredient, the name is sent as an input to the Instacart search bar
Step 3: Instacart looks for the best fit in all the stores nearby based on the ingredient name.
Step 4: Relying on the assumption, Selenium proceeds to select the 1st option that shows up.

The Instacart website did not showcase a particular schema for the price and the weight, However, scrolling down got the webpage in a format, from where we could get the specific details. This part took some time to figure out, without this step, scraping the prices would have been very cumbersome.

At this point, we have a dictionary of the ingredient names (Cluster name), which we used to map to the original ingredients. These results were then used to optimize cost while satisfying nutritional requirements.

### 3.2 The minimum and maximum nutritional component values - Scraping:
Below is a snippet of the maximum and minimum values of each of 72 nutritional components, scraped primarily from the National Institutes of Health website.[21]
(Provided in the Nutrient_values.xlsx file)

| Nutrient | Min_value | Max_value |
|---|---|---|
| 'total_fat (g)', | 44 | 78 |
| 'saturated_fat (g)', | 0 | 13 |
| 'cholesterol (mg)', | 0 | 200 |
| 'sodium (mg)', | 500 | 2000 |
| 'choline (mg)', | 10 | 450 |
| 'folate (mcg)', | 400 | 1000 |
| 'folic_acid (mcg)', | 400 | 1000 |
| 'niacin (mg)', | 16 | 35 |
| 'pantothenic_acid (mg)', | 5 | 1000 |
| 'riboflavin (mg)', | 1.3 | 1000 |
| 'thiamin (mg)', | 1.3 | 1000 |
| 'vitamin_a (IU)', | 3000 | 10000 |
| 'vitamin_a_rae (mcg)', | 900 | 3000 |
| 'carotene_alpha (mcg)', | 0 | 300000 |
| 'carotene_beta (mcg)', | 0 | 300000 |
| 'cryptoxanthin_beta (mcg)', | 0 | 300000 |
| 'lutein_zeaxanthin (mcg)', | 2 | 10 |
| 'vitamin_b12 (mg)', | 2.4 | 100 |
| 'vitamin_b6 (mg)', | 1.3 | 100 |
| 'vitamin_c (mg)', | 75 | 2000 |
| 'vitamin_d (IU)', | 600 | 4000 |
| 'vitamin_e (mg)', | 15 | 1000 |
| 'tocopherol_alpha (mg)', | 15 | 1000 |
| 'vitamin_k (mcg)', | 90 | 10000 |

Figure 5 : Nutritional values range as is required by the average male human body [21]

A major assumption here was that the min_value was taken as the minimum value for a female's diet and the max_value was taken as the maximum value for a male's diet; so we can get an encompassing range of values. These values were put together by brute force.

### 3.3 Linear Programming Problem:

The aim of the linear programming (LP) model in this context is to maximize the trade-off between the cost and nutritional content of a given set of products or ingredients. The aim function can be customized to either optimize the nutritional value within a certain budget or minimize the expenses while satisfying specific nutritional criteria.

### 3.3.1 Objective and Decision Variables
The objective function can be defined as follows:

- Maximizing Nutritional Value: Here, the objective is to maximize the total nutritional score, which is a weighted sum of various nutrients (like vitamins, minerals, proteins, etc.) present in the ingredients, subject to a constraint on the total cost not exceeding a predefined budget.

- Minimizing Cost: Conversely, the objective might be to minimize the total cost of the ingredients while ensuring that the nutritional content meets or exceeds a set of minimum required levels. This approach is particularly useful in scenarios where budget constraints are stringent, and the goal is to achieve the most nutritionally balanced diet within those financial limits.

- Weights: Allow fractional weights for ingredients by setting the variable type to 'GRB.CONTINUOUS'. This allows for more flexibility in the optimization problem.

### 3.3.2 Constraints:

- Maximum Nutrient Content: This represents the highest permissible level of each nutritional element in an ideal meal.

- Minimum Nutrient Content: This indicates the essential baseline amount of each nutritional element that must be met or exceeded in an ideal meal.

- Excess control Constraint: Some earlier models of our formulation gave exorbitant values for some particular ingredients. Notable example being 545 g out of ~800g of optimal solution being 'Cardoon'. Cardoons are dried leaves with high nutritional values. Although very healthy, we have to mitigate their amounts, since too much of one particular ingredient dominating the whole meal isn't realistically optimal.
  Thus, we add an excess control constraint akin to an 80-20 split. Where no particular ingredient can be more than 20% of the whole meal.

### 3.3.3 Final Formulation:

**Objective Function:**

Minimize the cost of the diet:

$$\text{Min.} \sum_{i \in ingredients} W_i \cdot P_i$$

Where:

- $W_i$ is the weight of ingredient i,
- $P_i$ is the price per unit weight of ingredient i.

**Constraints:**

1. **Nutrient Requirements**:
   For each nutrient k, ensure intake is within minimum $min_k$ and maximum $max_k$ bounds:

$$\sum_{i \in ingredients} W_i \cdot N_{i,k} \geq min_k \quad if \ k \in mini$$

$$\sum_{i \in ingredients} W_i \cdot N_{i,k} \leq max_k \quad if \ k \in maxi$$

Where $N_{i,k}$ is the nutrient k content in ingredient i; $maxi$ and $mini$ are collections of maximum and minimum nutrient requirements respectively.

2. **Diversity Constraint**:
Limit the contribution of any single ingredient to one-fifth of the total weight:

$$W_i \leq 1/5 \cdot \sum_j W_j \quad \forall \ i \ (ingredients)$$

3. **Non-Negativity**:
Weights of all ingredients must be non-negative:

$$W_i \geq 0 \quad \forall \ i \ (ingredients)$$

**Variables:**

- $W_i$: Continuous decision variable representing the weight of ingredient i.

# 4. Case Studies

We have identified five distinct scenarios in which our model proves invaluable by offering meticulously curated meal plans tailored to individuals with specific objectives. The outlined case studies encompass individuals with deficiencies in iron, and vitamin D. Additionally, two scenarios cater to those seeking a calorie deficit diet and a protein-rich diet, respectively. Furthermore, we have covered a hypothetical case of a maxed out diet.

**4.1 Optimal case (No Nutritional Deficiencies i.e. No Additional Constraints):**

The below values show the ingredients for the most cost-effective, nutritionally optimal diet that our model produces. This case, consider no nutritional deficiencies. We use the formulation, as it is, discussed in section 3.3.3. We can see at the bottom the cost of the diet(for the whole day) is $3.882, and the total weight of the diet is 783.27grams.

```
                                          Ingredient  Weight (g)
256                               Fish oil, cod liver    5.676798
285                                 Acerola juice, raw   0.592976
297                                 Onions, raw, welsh  33.332758
356                                   Fish oil, herring 17.873501
665                                           Yam, raw 156.653453
696                         Tofu, dried-frozen (koyadofu) 34.704419
883                          Seeds, low-fat, sesame flour 39.030415
2437               Nuts, blanched, hazelnuts or filberts   1.453651
2857       KELLOGG'S, Original 3-Bean Chips, BEANATURAL   9.573239
3598   Sea lion, meat with fat (Alaska Native), Steller  16.450028
3810     Cardoon, without salt, drained, boiled, cooked 156.653453
4291   Yam, without salt, or baked, drained, boiled, ...  35.291581
5058   Cereals, dry, 10 minute cooking, regular, CREA...  61.869963
5120   Cereals ready-to-eat, RALSTON Enriched Wheat B...   9.366579
5312   USDA Commodity Food, low saturated fat, vegeta...  10.399373
5334   Spaghetti, enriched (n x 6.25), dry, protein-f... 134.175633
5392   Wocas, yellow pond lily (Klamath), Oregon, dri...  12.748092
5929   Tofu, prepared with calcium sulfate, dried-fro...   9.67366
6220   Beef, raw, liver, variety meats and by-product...   0.54681
6525   Whale, skin and subcutaneous fat (muktuk) (Ala...  34.232059
8028   Margarine, CANOLA HARVEST Soft Spread (canola,...   2.968824

Total cost of the diet: $3.82212065153449

Total weight of the diet: 783.267265146858 grams
```

Figure 6 : Results of the Optimal Case

## 4.2 Vitamin D Deficiency

Vitamin D deficiency is a common global issue. About 1 billion people worldwide have vitamin D deficiency, while 50% of the population has vitamin D insufficiency; approximately 35% of adults in the United States have vitamin D deficiency[12]. Good food sources of Vitamin D include eggs, fish, and sardines. The best source of Vitamin D is sunlight. During the months of winter when the sun exposure is reduced, the population tends to experience a reduction in the general level of the vitamin levels [13]. A deficiency of Vitamin D can lead to loss of bone density, increased levels of stress, and in extreme cases it can lead to rickets [14].

An individual with a deficiency of Vitamin D requires 6000-10000 IU daily to overcome the deficiency [4].We created a curated ingredient list for a dietary plan for individuals with low levels of Vitamin D and below is the formulation and the results found by our model.

In the formulation, we make the following changes to the formulation discussed earlier in Section 3.3.3 - the nutrient requirements constraints are changed to the following constraints -

Minimum nutrient requirement:

- For Vitamin D ($k$ = "vitamin_d"):

$$\sum_{i \in ingredients} W_i \cdot N_{i,k} \geq 6000$$

- For other nutrients ($k \neq$ "vitamin_d"):

$$\sum_{i \in ingredients} W_i \cdot N_{i,k} \geq min_k$$

Maximum nutrient requirement:

- For Vitamin D ($k =$"vitamin_d"):

$$\sum_{i \in ingredients} W_i \cdot N_{i,k} \leq 10000$$

- For other nutrients ($k \neq$ "vitamin_d"):

$$\sum_{i \in ingredients} W_i \cdot N_{i,k} \leq max_k$$

Through this we get the following results:

```
                                              Ingredient  Weight (g)
180                                   Amaranth leaves, raw    7.871248
256                                      Fish oil, cod liver   9.056541
538                                           Nuts, almonds   0.531171
696                              Tofu, dried-frozen (koyadofu) 36.512829
883                             Seeds, low-fat, sesame flour  15.137087
1270                                   Nuts, dried, butternuts 30.877744
1291                                  Mushrooms, raw, maitake 207.865117
2027                                Fish, raw, Greenland, halibut 9.814195
2857              KELLOGG'S, Original 3-Bean Chips, BEANATURAL  95.653484
3162         Cereals, dry, plain, original, MALT-O-MEAL         19.149785
3372               Sea lion, liver (Alaska Native), Steller    10.966368
3598    Sea lion, meat with fat (Alaska Native), Steller       28.024233
3621         Sweeteners, packets, EQUAL, aspartame, tabletop   10.80851
4744   Sweetener, herbal extract powder from Stevia leaf       34.493254
4930    Cereals ready-to-eat, Warm Cinnamon, KASHI HEA...      11.56433
4947  Spaghetti, enriched (n x 6.25), cooked, protei...      207.865117
5058   Cereals, dry, 10 minute cooking, regular, CREA...      43.386985
5392   Wocas, yellow pond lily (Klamath), Oregon, dri...      21.397535
5617   Babyfood, organic, carrot and squash, apple, 2...       2.187064
5809   Pork, raw, frozen, ears, variety meats and by-...      17.733616
5929   Tofu, prepared with calcium sulfate, dried-fro...      10.564253
7436   Mushrooms, raw, exposed to ultraviolet light, ...     207.865117

Total cost of the diet: $12.30378978557047

Total weight of the diet: 1039.3255832078628 grams
```

Figure 7 : Results for the case of Vitamin D deficiency

**4.3 Iron deficiency**

Iron deficiency is a major issue in women and in particular pregnant women [7]. Iron deficiency anemia can lead to irregular heartbeat, enlarged hearts and heart failure. In pregnant women it can lead to premature births and low weight babies[15]. Good sources of iron are green leafy vegetables and red meat. For iron deficiency anemia pregnant women the recommended levels of iron intake is found to be at least 80 mg per day with an upper limit of the human body at 100mg per day[7]. Inputting these values into our model, we make the following changes to the formulation discussed earlier in Section 3.3.3 - the nutrient requirements constraints are changed to the following constraints -

Minimum nutrient requirement:

- For iron ($k$ = "iron"):

$$\sum_{i \in ingredients} W_i \cdot N_{i,k} \geq 80$$

- For other nutrients ($k \neq$ "iron"):

$$\sum_{i \in ingredients} W_i \cdot N_{i,k} \geq min_k$$

Maximum nutrient requirement:

- For iron ($k$ = "iron"):

$$\sum_{i \in ingredients} W_i \cdot N_{i,k} \leq 100$$

- For other nutrients ($k \neq$ "iron"):

$$\sum_{i \in ingredients} W_i \cdot N_{i,k} \leq max_k$$

Through this we get the following results:

```
                                Ingredient  Weight (g)
256                       Fish oil, cod liver    5.678974
285                         Acerola juice, raw    1.531905
297                         Onions, raw, welsh   24.347259
356                          Fish oil, herring    6.052562
665                                    Yam, raw  105.413056
696                 Tofu, dried-frozen (koyadofu)  41.309361
883                   Seeds, low-fat, sesame flour  34.240917
2085       Whale, raw (Alaska Native), meat, beluga 136.677879
2437         Nuts, blanched, hazelnuts or filberts   23.848323
2857     KELLOGG'S, Original 3-Bean Chips, BEANATURAL 29.953681
3162         Cereals, dry, plain, original, MALT-O-MEAL 2.016701
3598  Sea lion, meat with fat (Alaska Native), Steller  18.068058
3810    Cardoon, without salt, drained, boiled, cooked 163.252017
5058  Cereals, dry, 10 minute cooking, regular, CREA...  59.107738
5120  Cereals ready-to-eat, RALSTON Enriched Wheat B...   9.305103
5312  USDA Commodity Food, low saturated fat, vegeta...  20.201245
5334  Spaghetti, enriched (n x 6.25), dry, protein-f... 126.344594
5392  Wocas, yellow pond lily (Klamath), Oregon, dri...   2.570179
5929  Tofu, prepared with calcium sulfate, dried-fro...   6.340535

Total cost of the diet: $4.550389903280612

Total weight of the diet: 816.2600870215394 grams
```

Figure 8: Results for the case of Iron deficiency

**4.4 Calorie Deficit Diet**

Obesity is a huge concern in[16] today's world with about a third of the world's population being obese.

We have formulated a low-calorie diet with a consumption of 1000-1500 for those individuals who seek to follow a calorie-deficit diet. A deficit of 500 - 700 calories per day has been used for weight loss and is recommended by obesity guidelines and societies[17]. Inputting these values into our model, we make the following changes to the formulation discussed in Section 3.3.3 - the nutrient requirements constraints are changed to the following constraints -

Minimum nutrient requirement:

- For calories ($k$ = "calories (kcal)"):

$$\sum_{i \in ingredients} W_i \cdot N_{i,k} \geq 1000$$

- For other nutrients ($k \neq$ "calories (kcal)"):

$$\sum_{i \in ingredients} W_i \cdot N_{i,k} \geq min_k$$

Maximum nutrient requirement:

- For calories ($k$ = "calories (kcal)"):

$$\sum_{i \in ingredients} W_i \cdot N_{i,k} \leq 1500$$

- For other nutrients ($k \neq$ "calories (kcal)"):

$$\sum_{i \in ingredients} W_i \cdot N_{i,k} \leq max_k$$

Through this we get the following results:

```
                                        Ingredient  Weight (g)
180                             Amaranth leaves, raw    4.149511
256                                Fish oil, cod liver    5.80821
285                                  Acerola juice, raw   0.805296
297                                 Onions, raw, welsh  20.339133
665                                           Yam, raw 136.907198
696                       Tofu, dried-frozen (koyadofu)  25.702935
883                         Seeds, low-fat, sesame flour  26.830009
989                         Turnips, unprepared, frozen 160.264468
2437            Nuts, blanched, hazelnuts or filberts   50.305716
2857       KELLOGG'S, Original 3-Bean Chips, BEANATURAL   6.160871
3162        Cereals, dry, plain, original, MALT-O-MEAL   6.815593
3598  Sea lion, meat with fat (Alaska Native), Steller  23.656707
3810     Cardoon, without salt, drained, boiled, cooked 160.264468
3875      Tofu, prepared with calcium sulfate, firm, raw  31.917569
5120  Cereals ready-to-eat, RALSTON Enriched Wheat B...   5.236772
5334  Spaghetti, enriched (n x 6.25), dry, protein-f... 113.614812
5744  Cereals ready-to-eat, Honey Toasted Oat, KASHI...   1.669059
5929  Tofu, prepared with calcium sulfate, dried-fro...  14.236109
6157  Cereals ready-to-eat, KELLOGG'S ALL-BRAN COMPL...   0.80586
6295  Pork, braised, cooked, pancreas, variety meats...   4.961047
8491  Gelatin desserts, vitamin C, sodium, potassium...   0.870996

Total cost of the diet: $4.2584488171291

Total weight of the diet: 801.3223392730348 grams
```

Figure 9: Results for a calorie deficit diet

**4.5 High Protein Diet**

A protein-heavy diet is essential for various reasons. Proteins are vital for building and repairing tissues, supporting immune function, and maintaining muscle mass. They also play a crucial role in enzyme and hormone production, aiding in overall metabolic health. Incorporating ample protein ensures optimal body function, supports weight management, and promotes satiety, contributing to a balanced and a nutritious diet. Persons involved in heavy exercise and fitness regimes require even higher protein to maximize protein accretion and resistance exercise, which is about 1.6-2.2 grams of protein per kg of the body weight [18]. We input a minimum requirement of 128g and an upper limit of 184g of protein and make the following changes to the formulation discussed in Section 3.3.3 - the nutrient requirements constraints are changed to the following constraints -

Minimum nutrient requirement:

- For Protein ($k = $ "protein"):

$$\sum_{i \in ingredients} W_i \cdot N_{i,k} \geq 128$$

- For other nutrients ($k \neq $ "protein"):

$$\sum_{i \in ingredients} W_i \cdot N_{i,k} \geq min_k$$

Maximum nutrient requirement:

- For Protein ($k = $ "protein"):

$$\sum_{i \in ingredients} W_i \cdot N_{i,k} \leq 184.8$$

- For other nutrients ($k \neq $ "protein"):

$$\sum_{i \in ingredients} W_i \cdot N_{i,k} \leq max_k$$

Through this we get the following results:

```
                                          Ingredient  Weight (g)
256                                 Fish oil, cod liver    5.979223
285                                    Acerola juice, raw    3.043979
297                                    Onions, raw, welsh   22.649511
665                                              Yam, raw   85.08871
696                          Tofu, dried-frozen (koyadofu)  14.439823
883                           Seeds, low-fat, sesame flour  34.270467
2085          Whale, raw (Alaska Native), meat, beluga     19.499198
2857       KELLOGG'S, Original 3-Bean Chips, BEANATURAL    48.947594
3598   Sea lion, meat with fat (Alaska Native), Steller    5.159295
3810     Cardoon, without salt, drained, boiled, cooked  202.436824
4947  Spaghetti, enriched (n x 6.25), cooked, protei...  190.025636
5058     Cereals, dry, 10 minute cooking, regular, CREA...  59.265773
5120  Cereals ready-to-eat, RALSTON Enriched Wheat B...    0.602236
5312  USDA Commodity Food, low saturated fat, vegeta...   22.2083
5334  Spaghetti, enriched (n x 6.25), dry, protein-f...  202.436824
5392    Wocas, yellow pond lily (Klamath), Oregon, dri...   33.630926
5929  Tofu, prepared with calcium sulfate, dried-fro...    7.653978
6295  Pork, braised, cooked, pancreas, variety meats...   10.371225
6525  Whale, skin and subcutaneous fat (muktuk) (Ala...   44.474599

Total cost of the diet: $4.338695988601508

Total weight of the diet: 1012.1841186211607 grams
```

Figure 10: Results for a high protein diet

### 4.6 God's Diet (Maxed-Out Diet):

We define here a hypothetical case, where we are primarily concerned with absolutely maxing out on all the nutritional components' value of our ingredients(still within bounds), while only minimizing the price of the meal as a secondary objective. For instance, the protein requirement is between 60 gm-200 gm for a day, and this particular instance of the model tried to push the protein weight in the diet to as close to 200gm as possible without a lot of care for the price.

Note: By altering the diet constraints in this case, individuals can tailor their meal plans to prioritize optimal nutrition while also having the option to apply specific cost and weight limitations, ensuring a diet that balances health goals with budget and portion considerations.

### 4.6.1 Gurobi's Multi-Objective Environment:
We make use of Gurobi's multi-objective environment to set multiple objectives:
1. Objective 1(Obj1): Maximize the nutritional components' values
2. Objective 2(Obj2): Minimize the price.

And we make use of the corollary for maxing out the nutritional components' value:

$$\text{Min. } f(x) \equiv \text{Max. } -f(x)$$

.(since the ModelSense, can only take either GRB.MINIMIZE or GRB.MAXIMIZE)

To incorporate all of this, we made the following changes to our formulation:

- **Primary Objective** (i.e. Obj1: Maximization of Nutrient Value):

Maximize the total nutrient value across all ingredients and nutrients:

$$Maximize: \sum_{i \in ingredients} \sum_{k \in keys} W_i \cdot N_{i,k}$$

This is implemented in the model as:

$$Minimize: -\sum_{i \in ingredients} \sum_{k \in keys} W_i \cdot N_{i,k}$$

- **Secondary Objective** (i.e. Obj2: Minimization of Price):

Minimize the total cost of the diet:

$$Minimize: \sum_{i \in ingredients} W_i \cdot P_i$$

Further, we can assign the priority to the two objectives in our code as:

Obj1 gets priority = 1, i.e. primary.
```
model.setObjectiveN(obj1, 0, 1)
```

And, Obj2 gets priority = 0, i.e. secondary.
```
model.setObjectiveN(obj2, 1, 0)
```

Through this we get the following results:
```
                                          Ingredient    Weight (g)
67                              Oil, soybean lecithin      7.762224
180                                Amaranth leaves, raw    1.286512
256                                    Fish oil, cod liver 1.221551
450                                         Oil, wheat germ    4.36787
711                        Pumpkin, without salt, canned    52.305005
1378              Soy protein isolate, potassium type     19.088659
1618                 Gelatins, unsweetened, dry powder     40.387709
1761                    Cereals, dry, unenriched, farina    69.155478
2057          SMART SOUP, Vietnamese Carrot Lemongrass    154.707742
2140              KASHI, unprepared, 7 Whole Grain, Pilaf   86.243194
2601                 Beverages, sugar free, Energy Drink   402.474342
2775       KASHI, Frozen Entree, Chicken Pasta Pomodoro   113.972411
2828       Beverages, decaffeinated, brewed, green, tea 87684.987617
3372             Sea lion, liver (Alaska Native), Steller   42.70884
3598    Sea lion, meat with fat (Alaska Native), Steller   58.387807
4285   Bamboo shoots, without salt, drained, boiled, ...  460.13465
4480   Beverages, Glacial Natural spring water, ICELA... 87956.50965
4583   Beverages, brewed, other than chamomile, herb,...  20239.6774
4668   Beverages, DANNON, non-carbonated, bottled, water 87956.50965
4671   Candies, dietetic or low calorie (sorbitol), hard  151.289218
4744   Sweetener, herbal extract powder from Stevia leaf  596.310724
4899   Babyfood, without added fluoride., GERBER, bot... 87956.50965
5617   Babyfood, organic, carrot and squash, apple, 2...   58.823529
5767   Beverages, fortified with vitamin C, Apple jui...  406.436241
6209   Alcoholic beverage, all (gin, rum, vodka, whis...   70.528967
6272   Beverages, diet, peach, ready to drink, black ...  3390.102631
7403   Beverages, non-carbonated, bottled, water, DAS... 61583.347213
7436   Mushrooms, raw, exposed to ultraviolet light, ...  161.693379
7774   Oil, woks and light frying, principal uses sal...   41.822239
7839   Beverages, fortified with vitamin C, unsweeten...   13.796147

Total cost of the diet: $87927.79144622775

Total weight of the diet: 439782.5482489649 grams
```

Figure 11: The results obtained from utilising Gurobi multi-objective environment

*Note: If we switch the priorities of the Obj1 and Obj2 i.e. put minimizing the price as the primary objective, we get the results same as the aforementioned "Optimal case" referenced in 4.1.*

## 5. Results

Based on the various cases shown though the model we formulate some inferences from the results which are discussed and compared below.

## 5.1 Comparison based on weight

Based on a Gurobi optimisation result, the scatter plot "Weight of meal (g) vs. Dietary preference" shows the meal weights for different dietary needs. There seems to be a relationship between dietary requirements and meal weight because the heaviest meal is for "Vitamin D Deficiency" at 1039.325 grams, while the "Optimal Case" meal weighs 783.267 grams. The "Protein High" and "Iron Deficiency Anaemia" diets weigh 816.26 and 1012.184 grams, respectively, whereas the "Calorie Deficit" diet weighs 801.322 grams. The graph shows that meal content and portion size are highly influenced by nutritional preferences.

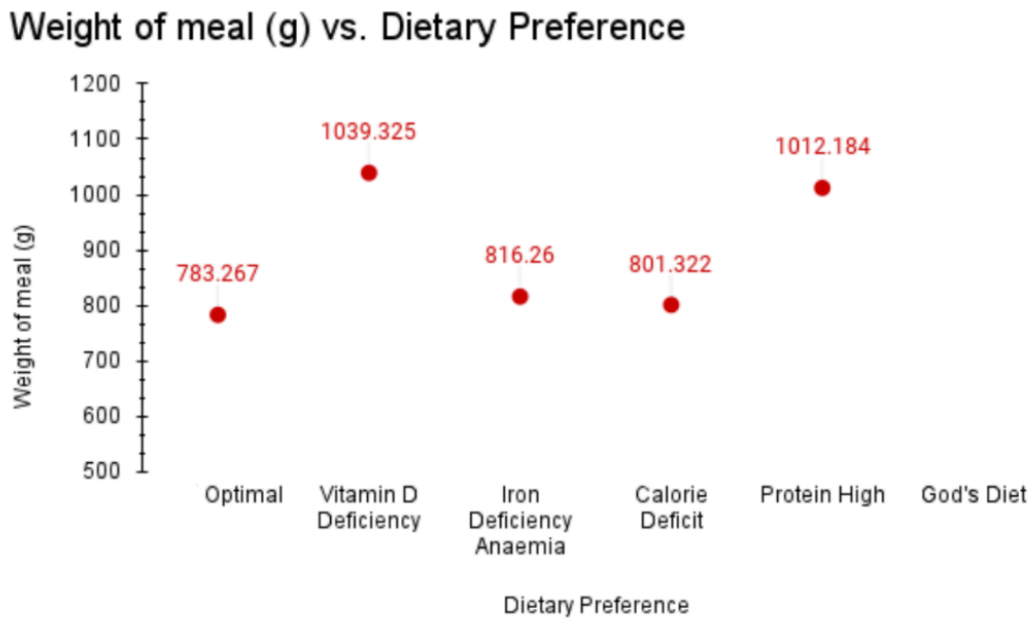

Figure 12: Scatter plot to represent the comparison based on weight

## 5.2 Comparison based on cost

A comparison of meal costs based on a Gurobi optimisation model is shown in the bar graph below, which shows a significant price difference for diets that address vitamin D insufficiency. This surge suggests that compared to other dietary situations, obtaining the recommended amount of Vitamin D may require more expensive or specialized products. The cost-effectiveness of the model for other dietary preferences seems to be fairly balanced, underscoring the difficulty in treating vitamin D insufficiency in diet planning from an economic standpoint. The average cost of a nutritious meal for an adult no matter what their dietary requirements is found to be $5.84. This shows the minimum possible cost for a whole day's nutrition for an average adult.

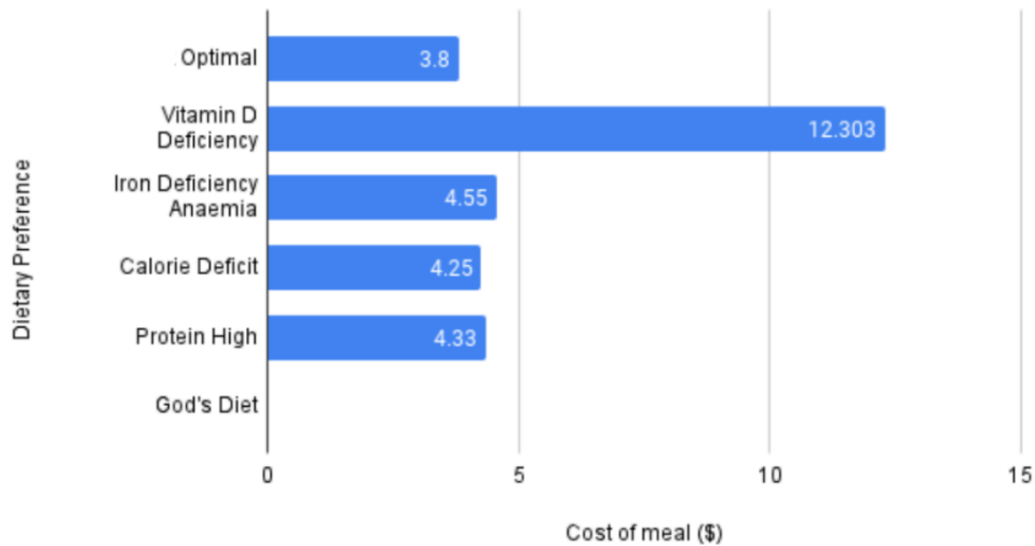

Figure 13 : Graphical representation for comparison of cases based on cost

**5.3 Analysis of a normal (the optimal) diet**

Below we show a pie chart representation of the ingredients suggested by our model in a normal case scenario for a healthy adult with no additional requirements. We see that three ingredients form a majority of the diet while the rest are occupied in smaller proportions by other food items. This shows a balanced diet where micros and the macros of one's diet are handled. It also nearly follows the rules of thirds of diet[23].

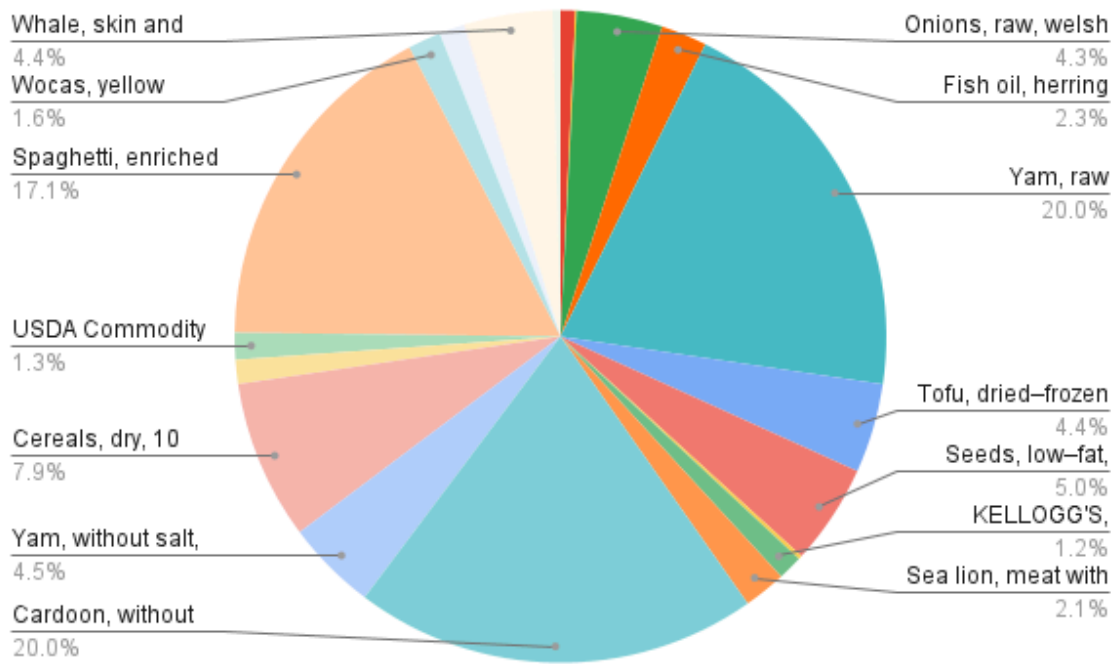

Figure 14: Pie chart analysis of our optimal diet

Following is the optimal list of ingredients for a diet with the optimal nutritional value and minimum price for a normal diet of an adult with no additional nutritional requirements.

**5.4 Time Complexity Analysis**

Dual-Simplex iterations vs time graph vs the cases we used. We formulate a comparison between the time taken for each case and the number of iterations our model takes to provide us the result. We calculate the correlation coefficient between the two values and it comes to -0.236. This indicates a weak relation between the two and not a direct one. It means that as time increases there is a slight tendency for the number of iterations to increase but it is not a strong relationship. Further, it needs to be pointed out that the most general case (that being corresponding to god's diet) takes about 10 minutes to run, and all other takes about an average of 300 seconds, indicating the celerity of the model.

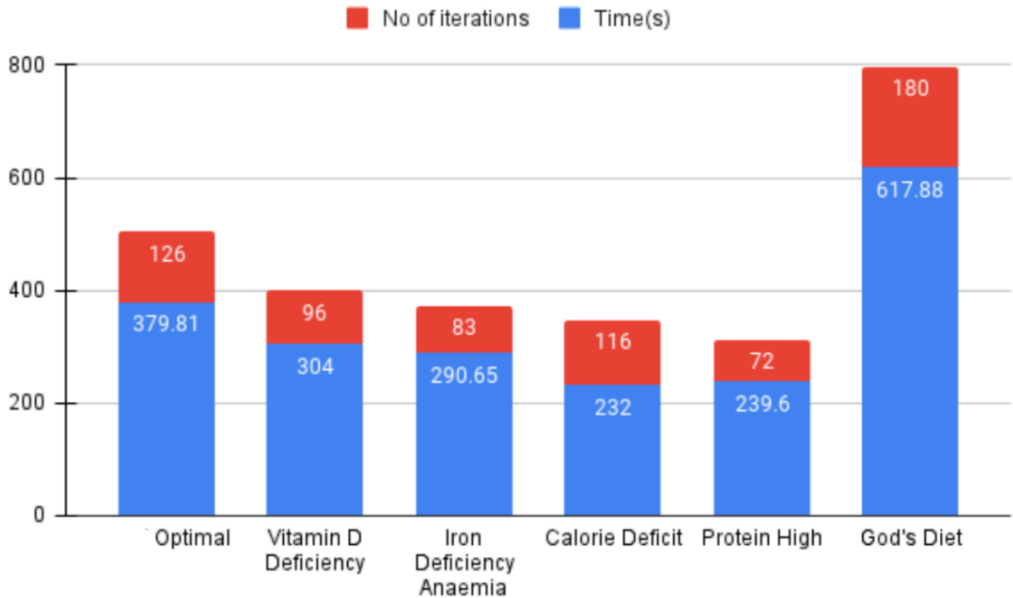

Figure 15 : Visualization of time complexity analysis for all the discussed cases

## 6. Conclusions

The model and its formulation effectively leverages linear programming, making use of Gurobi optimization tools and libraries across a meticulously compiled dataset. This integration results in a valuable solution that provides the necessary ingredients for creating the optimal nutritional diet at the lowest cost.

This model finds utility in its potential to save time and streamline decision-making in the creation of health-conscious meals that minimize the price.

The model showcases adaptability, allowing for customization to several nutritional deficiencies by introducing just a few new constraints. While we've only discussed its potential for customization for addressing these deficiency cases in this report, it's important to note that the model can be similarly tailored to individual preferences, ingredient availability, relaxed price constraints, and more.

In real-life scenarios, it serves as an excellent means for establishing a lower limit for an optimal diet within the constraints of nutritional considerations, including special dietary needs. For instance, consider an individual suffering from Iron Deficiency Anaemia (IDA), as explored in the case studies.

This model proves invaluable in determining the minimum monthly meal cost required to address their Iron deficiency. To illustrate, 30 days at $4.55 each amounts to $136.5/month.

## 7. Future Work

We can investigate methods for collecting detailed data on individual taste satisfaction levels. This would involve gathering subjective taste preferences from a variety of participants, which can then be integrated into our diet planning model (such as a particular demographic or society). This data will enhance the personalization of the diet plan, ensuring it aligns not only with nutritional and cost goals but also with individual taste preferences (different cuisines).

Integrating our personalized diet planning tool into an inventory system would significantly streamline decision-making for nutritionists, doctors, and individuals alike. This system would offer real-time tracking of ingredient availability, align dietary recommendations with current stock, and adapt to changing nutritional needs and preferences. Such a tool would not only facilitate efficient meal planning and shopping but also empower users with actionable insights for maintaining a balanced diet, making nutritional management more accessible and effective for everyone.

While the outputs of these results may not always align with a conventional recipe for preparing a dish/dishes, the model developed remains valuable especially when augmented with an additional layer of Natural Language Processing.

These enhancements can be implemented by loosening the constraints of the formulated LP, thereby accommodating a broader range of ingredients within the optimal solution. A diverse array of dishes data can be systematically scraped from online sources. Subsequently, employing LLM(s) we can form suggestions that map all the ingredients from our optimal solution to feasible dishes that can be prepared rather easily. This mapping process can be refined to minimally impact the optimal price of the diet or its nutritional value. A framework and methodology akin to the aforementioned is presented in the work - Willis, A. (2017) [22].

# Bibliography

[1] Wanjek, C. (2005). Food at Work: Workplace solutions for malnutrition, obesity and chronic diseases.

[2] Klímová, B., Dziuba, S., & Cierniak–Emerych, A. (2020). The Effect of Healthy Diet on Cognitive Performance Among Healthy Seniors – A Mini Review. *Frontiers in Human Neuroscience*, *14*. https://doi.org/10.3389/fnhum.2020.00325

[3] Cena, H., & Calder, P. C. (2020). Defining a Healthy Diet: Evidence for the role of contemporary dietary patterns in Health and disease. *Nutrients*, *12*(2), 334. https://doi.org/10.3390/nu12020334

[4] Sizar, Omeed, et al. "Vitamin D Deficiency." *PubMed*, StatPearls Publishing, 17 July 2023, www.ncbi.nlm.nih.gov/books/NBK532266/.

[5] Türkmenoğlu, Cumali, et al. "Recommending Healthy Meal Plans by Optimising Nature-Inspired Many-Objective Diet Problem." *Health Informatics Journal*, vol. 27, no. 1, Jan. 2021, p. 146045822097671, https://doi.org/10.1177/1460458220976719

[6] Pochmann, Vitor O, and Fernando J Von. "Multi-Objective Bilevel Recommender System for Food Diets." *2022 IEEE Congress on Evolutionary Computation (CEC)*, 18 July 2022, https://doi.org/10.1109/cec55065.2022.9870408

[7] Kumar, Aditi, et al. "Iron Deficiency Anaemia: Pathophysiology, Assessment, Practical Management." *BMJ Open Gastroenterology*, vol. 9, no. 1, 1 Jan. 2022, p. e000759, bmjopengastro.bmj.com/content/9/1/e000759, https://doi.org/10.1136/bmjgast-2021-000759 .

[8] Kim, Ju Young. "Optimal Diet Strategies for Weight Loss and Weight Loss Maintenance." *Journal of Obesity & Metabolic Syndrome*, vol. 30, no. 1, 27 Oct. 2020, www.ncbi.nlm.nih.gov/pmc/articles/PMC8017325/, https://doi.org/10.7570/jomes20065.

[9] *Nutritional values for common foods and products*. (2019, June 18). Kaggle. https://www.kaggle.com/datasets/trolukovich/nutritional-values-for-common-foods-and-products

[10] Rowe, Sam, and Anitra C. Carr. "Global Vitamin c Status and Prevalence of Deficiency: A Cause for Concern?" *Nutrients*, vol. 12, no. 7, 6 July 2020, p. 2008, www.ncbi.nlm.nih.gov/pmc/articles/PMC7400810/, https://doi.org/10.3390/nu12072008.

[11] Léger, Daniel. "Scurvy: Reemergence of Nutritional Deficiencies." *Canadian Family Physician Medecin de Famille Canadien*, vol. 54, no. 10, 2008, pp. 1403–6, www.ncbi.nlm.nih.gov/pmc/articles/PMC2567249/#:~:text=Treatment%20and%20prognosis.

[12] Cleveland Clinic. "Vitamin D Deficiency: Causes, Symptoms & Treatment." *Cleveland Clinic*, 8 Feb. 2022, my.clevelandclinic.org/health/diseases/15050-vitamin-d-vitamin-d-deficiency

[13] Hansen, Anita L., et al. "Vitamin D Supplementation during Winter: Effects on Stress Resilience in a Randomized Control Trial." *Nutrients*, vol. 12, no. 11, 24 Oct. 2020, p. 3258, https://doi.org/10.3390/nu12113258.

[14]"Vitamin D Deficiency." *Medlineplus.gov*, medlineplus.gov/vitaminddeficiency.html#:~:text=What%20problems%20does%20vitamin%20D

[15] Mayo Clinic. "Iron Deficiency Anemia - Symptoms and Causes." *Mayo Clinic*, Mayo Clinic, 4 Jan. 2022, www.mayoclinic.org/diseases-conditions/iron-deficiency-anemia/symptoms-causes/syc-20355034.

[16]Koliaki, Chrysi, et al. "Defining the Optimal Dietary Approach for Safe, Effective and Sustainable Weight Loss in Overweight and Obese Adults." *Healthcare*, vol. 6, no. 3, 28 June 2018, p. 73, www.ncbi.nlm.nih.gov/pmc/articles/PMC6163457/, https://doi.org/10.3390/healthcare6030073.

[17]Kim, Ju Young. "Optimal Diet Strategies for Weight Loss and Weight Loss Maintenance." *Journal of Obesity & Metabolic Syndrome*, vol. 30, no. 1, 27 Oct. 2020, www.ncbi.nlm.nih.gov/pmc/articles/PMC8017325/, https://doi.org/10.7570/jomes20065.

[18]Stokes, Tanner, et al. "Recent Perspectives Regarding the Role of Dietary Protein for the Promotion of Muscle Hypertrophy with Resistance Exercise Training." *Nutrients*, vol. 10, no. 2, 7 Feb. 2018, p. 180, www.mdpi.com/2072-6643/10/2/180/pdf, https://doi.org/10.3390/nu10020180.

[19]Kiani, Aysha Karim. "Main Nutritional Deficiencies." *Journal of Preventive Medicine and Hygiene*, vol. 63, no. 2 Suppl 3, 1 June 2022, pp. E93–E101, pubmed.ncbi.nlm.nih.gov/36479498/, https://doi.org/10.15167/2421-4248/jpmh2022.63.2S3.2752.

[20]*WebDriver*. (n.d.). Selenium. https://www.selenium.dev/documentation/webdriver/

[21]*Dietary supplement fact sheets*. (n.d.). https://ods.od.nih.gov/factsheets/list-all/

[22]Willis, A. (2017). *Forage: Optimizing food use with machine learning generated recipes*. https://www.semanticscholar.org/paper/Forage%3A-Optimizing-Food-Use-With-Machine-Learning-Willis-Lin/02478e392bce7167ea70f7ddbabc5ec06c090026

[23]"Easy Meal Planning: The Rule of Thirds to Build Your Body." *LeanBody.com*, leanbody.com/blogs/protein-corner/easy-meal-planning-the-rule-of-thirds-to-build-your-body. Accessed 12 Dec. 2023.